\documentclass[12pt,a4paper]{article}
\usepackage{amssymb}

\usepackage{amsmath}
\usepackage{color}
\usepackage[colorlinks,linkcolor=blue,citecolor=red]{hyperref}

\newtheorem{theorem}{Theorem}

\begin{document}

\title{Integrable  equations of the dispersionless Hirota type and hypersurfaces in the Lagrangian Grassmannian}

\author{E.V. Ferapontov, \  L. Hadjikos and K.R. Khusnutdinova}
     \date{}
     \maketitle
     \vspace{-5mm}
\begin{center}
Department of Mathematical Sciences \\ Loughborough University \\
Loughborough, Leicestershire LE11 3TU \\ United Kingdom \\[2ex]
e-mails: \\[1ex] \texttt{E.V.Ferapontov@lboro.ac.uk}\\
\texttt{L.Hadjikos@lboro.ac.uk}\\
\texttt{K.Khusnutdinova@lboro.ac.uk}
\end{center}

\vspace{1cm}

\begin{abstract}
We investigate integrable second order equations
of the  form
$$
F(u_{xx}, u_{xy}, u_{yy}, u_{xt}, u_{yt}, u_{tt})=0.
$$
Familiar examples include the Boyer-Finley equation
$u_{xx}+u_{yy}=e^{u_{tt}}$, the potential form of the dispersionless Kadomtsev-Petviashvili  (dKP) equation $u_{xt}-\frac{1}{2}u_{xx}^2=u_{yy}$, the dispersionless Hirota equation $(\alpha -\beta)e^{u_{xy}}+
(\beta -\gamma)e^{u_{yt}}+(\gamma -\alpha)e^{u_{tx}}=0$, etc.
The integrability is understood as  the existence of infinitely many hydrodynamic reductions.
We demonstrate that the natural equivalence group of the problem is isomorphic to $Sp(6)$, revealing a remarkable correspondence between differential equations of the above type and hypersurfaces of the Lagrangian Grassmannian. We prove that the moduli space of integrable equations of the dispersionless Hirota type is $21$-dimensional, and the action of the equivalence group $Sp(6)$ on the moduli space has an open orbit.

\bigskip


\bigskip

\end{abstract}

\newpage

\section{Introduction}

We investigate a general class of three-dimensional second order equations of the form
\begin{equation}
F(u_{xx}, u_{xy}, u_{yy}, u_{xt}, u_{yt}, u_{tt})=0
\label{main}
\end{equation}
where $u=u(x, y, t)$ is a function of three independent variables. Equations of this type arise in a wide range of applications including non-linear physics, general relativity, differential geometry, integrable systems and complex analysis. For instance, the dKP equation, 
$$
u_{xt}-\frac{1}{2}u_{xx}^2=u_{yy},
$$ 
also known as the Khokhlov-Zabolotskaya equation, arises in  non-linear acoustics \cite{KZ}, as well as in the theory of Einstein-Weyl structures \cite{D}. The Boyer-Finley equation,
$$
u_{xx}+u_{yy}=e^{u_{tt}},
$$
which  is descriptive of a class of self-dual $4$-manifolds,  has been extensively discussed  in the context of general relativity \cite{BF}. The equations
$$
Hess~u=1 ~~~ {\rm  and } ~~~ Hess~u=\triangle u,
$$
where $Hess$ is the determinant of the Hessian matrix of $u$, and $\triangle$ is the Laplacian, appear in differential geometry in the theory of  affine spheres and special Lagrangian $3$-folds, respectively \cite{Calabi, Joyce}.  A  subclass of equations of the form (\ref{main}),
$$
u_{tt}=f(u_{xx}, u _{xt}, u _{xy}),
$$
was discussed recently in \cite{MaksEgor} in connection with  hydrodynamic chains satisfying the so-called Egorov  property. Equations of the form (\ref{main}) typically arise as the Hirota-type relations for  various  $(2+1)$-dimensional  dispersionless hierarchies. For instance, equations of the  dispersionless  Toda hierarchy \cite{Takasaki, Kodama, Zab},
$$
\begin{array}{c}
(\lambda -\mu)e^{D_{\lambda}D_{\mu}u}=\lambda e^{-\partial_sD_{\lambda}u}-\mu e^{-\partial_sD_{\mu}u}, \\
\ \\
(\lambda -\mu)e^{{\bar D}_{\lambda}{\bar D}_{\mu}u}=\lambda e^{\partial_s{\bar D}_{\mu}u}-\mu e^{\partial_s{\bar D}_{\lambda}u}, \\
\ \\
e^{D_{\lambda}{\bar D}_{\mu}u}=1-\frac{\mu}{\lambda}e^{(\partial_s^2+\partial_sD_{\lambda}-\partial_s{\bar D}_{\mu})u},
\end{array}
$$
 take the form (\ref{main}) after one replaces the `vertex' operators
$D_{\lambda}, \ D_{\mu}, \ {\bar D}_{\lambda}, \ {\bar D}_{\mu}$ by ordinary partial derivatives. Equations of the dispersionless Toda hierarchy play an  important role in complex analysis underlying the integrable structure of the Dirichlet  problem for simply-connected domains  \cite{Wiegmann,  Marshakov}.  A generalization of this construction to multiply-connected domains \cite{Krich+Zabr} leads to a remarkable hierarchy of Hirota-like equations which constitute the so-called universal Whitham hierarchy \cite{Kr}.
Further examples of  Hirota-type relations arise in the theory of the  associativity (WDVV) equations \cite{B, Braden, H}, e.g.,
$$
e^{u_{xx}+u_{xy}+u_{xt}-u_{yt}}-e^{u_{xy}+u_{yy}+u_{yt}-u_{xt}}+e^{u_{xt}+u_{yt}+u_{tt}-u_{xy}}=0,
$$
$$
e^{u_{xx}+2u_{xt}+u_{xy}}-e^{u_{tt}+2u_{xt}+u_{yt}}=1,
$$
$$
e^{u_{tt}-u_{xx}}= \sinh u_{xy} /  \sinh u_{yt},
$$
$$
\coth u_{xy}=\coth u_{xt} \coth u_{yt},
$$
etc. These equations arise as certain differential constraints which should be considered along with a full set of  the WDVV equations.

Equations of the above type have been approached by a whole variety of modern techniques including symmetry analysis, differential-geometric and algebro-geometric methods, dispersionless $\bar \partial$-dressing, factorization techniques, Virasoro constraints, hydrodynamic reductions, etc. However, until recently there was no `intrinsic'  approach which would explain the integrability of these and other examples.  Moreover, there was no satisfactory definition of the integrability  which would (a) be algorithmically verifiable, (b) allow classification results and (c) provide a scheme for the construction of exact solutions. We emphasize that equations of the form (\ref{main}) are not amenable to the inverse scattering transform, and require an alternative approach. 
Such approach, based on the method of hydrodynamic reductions and, primarily,  the work \cite{GibTsar}, see also \cite{Kodama, F, Fer22, M}, etc, was proposed in \cite{Fer1}. It was suggested to define the integrability of a multi-dimensional dispersionless system by requiring the existence of `sufficiently many' hydrodynamic reductions which provide multi-phase solutions known as non-linear interactions of planar simple waves. Technically, one `decouples' a three-dimensional PDE (\ref{main}) into a pair of commuting $n$-component $(1+1)$-dimensional systems of hydrodynamic type,
\begin{equation}
R_{t}^{i}=\lambda ^{i}(R)\ R_{x}^{i}, \ \ \ \ \  R_{y}^{i}=\mu^{i}(R)\ R_{x}^{i},
\label{R}
\end{equation}
where the characteristic speeds $\lambda^i$ and $\mu^i$ satisfy the commutativity conditions 
\begin{equation}
\frac{\partial _{j}\lambda ^{i}}{\lambda ^{j}-\lambda ^{i}}=\frac{\partial
_{j}\mu ^{i}}{\mu ^{j}-\mu ^{i}}, ~~~ i\ne j,  \label{Com}
\end{equation}
$\partial_j=\partial_{R^j}, $ see \cite{Tsar}. 
The best way to illustrate the method of hydrodynamic reductions is to discuss an example.

{\bf Example. } Let us consider the dKP equation, $u_{xt}-\frac{1}{2}u_{xx}^2=u_{yy}, $
 and introduce the notation 
$
u_{xx}=a,\ u_{xy}=b, \ u_{xt}=p, \ u_{yy}=p-\frac{1}{2}a^2;
$
this results in the equivalent quasilinear representation of the dKP equation,
\begin{equation}
a_y=b_x, ~~ a_t=p_x,  ~~ b_t=p_y, ~~ b_y=(p-\frac{1}{2}a^2)_x.
\label{eq}
\end{equation}
Let us seek multi-phase solutions in the form
$a=a(R^1, ..., R^n), \ b=b(R^1, ..., R^n), \  p=p(R^1, ..., R^n)$
where the `phases'  $R^i(x, y, t)$ satisfy the equations (\ref{R}). The substitution of this ansatz into (\ref{eq}) implies the relations
$$
\partial_i b=\mu^i \partial_i a, ~~~ \partial_i p=\lambda^i \partial_i a, ~~~ \lambda^i=a+(\mu^i)^2.
$$
Calculating the compatibility conditions
$\partial_i\partial_jb=\partial_j\partial_ib$, $\partial_i\partial_jp=\partial_j\partial_ip$,  and substituting 
$ \lambda^i=a+(\mu^i)^2$ into the commutativity conditions
 (\ref{Com}), one obtains the  
Gibbons-Tsarev system for $a(R)$ and $\mu^i(R)$,
$$
\partial_j\mu^i=\frac{\partial_j a}{\mu^j-\mu^i}, ~~~
\partial_i\partial_ja=2\frac{\partial_ia\partial_ja}{(\mu^j-\mu^i)^2},
$$
$i\ne j$, which was first derived in \cite{GibTsar} in the
theory of hydrodynamic reductions of
Benney's moment equations. It is  remarkable  that the Gibbons-Tsarev system is in involution, and its general solution depends, modulo reparametrizations $R^i\to \varphi^i(R^i)$, on $n$ arbitrary functions of one variable. Thus, the dKP equation possesses infinitely many $n$-component reductions parametrized by $n$ arbitrary functions of one variable.  We point out that the compatibility conditions 
$\partial_k\partial_j\mu^i =\partial_j\partial_k\mu^i$ and  $\partial_i\partial_j\partial_k a=\partial_i\partial_k\partial_j a$     involve {\it triples} of indices $i\ne j \ne k$ only. Thus, for $n=2$ the Gibbons-Tsarev system  is automatically consistent, while its consistency   for $n=3$ implies the consistency for arbitrary $n$. 
Based on this example, we give the following

\noindent {\bf Definition} {\it An equation of the form (\ref{main}) is said to be integrable if, for any $n$, it possesses infinitely many $n$-component hydrodynamic reductions parametrized by $n$ arbitrary functions of one variable. }

We have verified that all examples presented above are indeed integrable in this sense, with the exception of the equations $Hess~u=1$ and $Hess~u=\triangle u$, which do not pass the test (see Sect. 3.5). In Sect.  2 we derive the  integrability conditions  (as a system of third order differential relations for the function $F$ in (\ref{main})), and prove our first main result (Theorem 1 of Sect. 2):

\begin{itemize}

\item {\bf The moduli space of integrable equations of the dispersionless Hirota type is $21$-dimensional.}

\end{itemize}

\noindent  The class of equations (\ref{main})  is form-invariant under the action of the contact group $Sp(6)$ generated by linear symplectic transformations of the variables $x, y, t, u_x, u_y, u_t$. These transformations  map integrable equations to integrable. Taking into account  that dim $Sp(6)=21$, our second main result (Theorem 5 of Sect. 6) reads as follows:

\begin{itemize}

\item  {\bf The action of the equivalence group $Sp(6)$ on the moduli space of integrable equations  has an open orbit.}

\end{itemize}

\noindent This result is, in a sense, surprising: it establishes the existence of a unique `master-equation'  which generates all other integrable examples via various (singular) limits. The existence of an open orbit does not mean, of course, that all integrable equations of the form (\ref{main}) are $Sp(6)$-equivalent. The structure of the orbit space  of this action is quite complicated, for instance, all examples listed above belong to various singular orbits of lower dimension, and are not $Sp(6)$-equivalent. 

From differential-geometric point of view, the equation (\ref{main}) can be viewed as defining a $5$-dimensional hypersurface $M^5$ in the Lagrangian Grassmannian $\Lambda$ which can (locally) be identified with the space of $3\times 3$ symmetric matrices  $u_{ij}$. In Sect. 4 we demonstrate that each hypersurface $M^5\subset \Lambda$ inherits a special structure, namely, a field of rational normal curves $\gamma $ of degree four specified in the tangent bundle $TM^5$. Bisecant surfaces $\Sigma^2\subset M^5$ are defined as two-dimensional submanifolds  whose projectivised tangent spaces are bisecant lines of $\gamma$. 
Similarly, trisecant  $3$-folds are three-dimensional submanifolds $\Sigma^3\subset M^5$ whose projectivised tangent spaces are trisecant planes of $\gamma$.  Our third main result is the following
geometric characterization of integrable equations (Theorem 3 of Sect. 4):

\begin{itemize}

\item {\bf Bisecant surfaces and trisecant $3$-folds are geometric images of the two- and three-component hydrodynamic reductions, respectively. In particular, the equation  is integrable if and only if the corresponding hypersurface $M^5$ possesses infinitely many  trisecant $3$-folds.}

\end{itemize}

The paper is organized as follows.

In Sect. 2 we derive the integrability conditions  by  transforming (\ref{main}) to a quasilinear form, and  calculating the compatibility conditions of the corresponding generalized Gibbons-Tsarev system which governs $n$-component hydrodynamic reductions. The result of this calculation is a complicated set of the integrability conditions, which are in involution. This establishes the $21$-dimensionality of the moduli space.

Further examples of integrable equations and partial classification results are provided in Sect 3, based on the integrability conditions obtained in Sect. 2. In particular, we analyze  the integrability of the symplectic Monge-Amp\'ere equations.

 Differential-geometric aspects of integrable equations of the form (\ref{main}) are discussed in Sect. 4.
 We begin by introducing a flat generalized conformal structure  on the  Lagrangian Grassmannian $\Lambda$, whose   conformal automorphism group is isomorphic to $Sp(6)$. Given a hypersurface $M^5\subset \Lambda$, we describe  geometric structures  which are naturally induced on $M^5$, and provide a geometric characterization of the integrability conditions. 
 
 In Sect. 5 we briefly discuss  multi-dimensional equations of the form (\ref{main}).
 
 The action of the equivalence group $Sp(6)$ on the moduli space of integrable equations  is discussed in Sect. 6.

\section{Derivation of the integrability conditions: proof of Theorem 1}

In this section we derive the integrability conditions which follow from the requirement of the existence of $n$-component reductions, and prove the following

\begin{theorem} The moduli space of integrable equations of the dispersionless Hirota type is $21$-dimensional. 
\end{theorem}

\centerline{\bf Proof:}

To apply the method of hydrodynamic reductions we first rewrite the  equation (\ref{main}) in the evolutionary form,
\begin{equation}
u_{tt}=f(u_{xx}, u_{xy}, u_{yy}, u_{xt}, u_{yt}),
\label{evol}
\end{equation}
and introduce the notation
$$
u_{xx}=a,\ u_{xy}=b, \ u_{yy}=c,\ u_{xt}=p,\ u_{yt}=q, \ u_{tt}= f(a,b, c, p, q).
$$
This provides an equivalent  quasilinear representation of (\ref{evol}),
\begin{equation}
\begin{array}{c}
a_y=b_x, ~~ a_t=p_x, ~~ b_y=c_x, ~~ b_t=p_y=q_x, ~~ c_t=q_y, \\
\ \\
p_y=q_x, ~~ p_t=f(a,b, c, p, q)_x, ~~ q_t=f(a,b, c, p, q)_y.
\end{array}
\label{quasi}
\end{equation}
Looking for solutions in the form $a=a(R^{1},...,R^{n}),\ b=b(R^{1},...,R^{n}),\ c=c(R^{1},...,R^{n}),\ p=p(R^{1},...,R^{n}), \ q=q(R^{1},...,R^{n})$, where the Riemann invariants $R^i$ satisfy the equations (\ref{R}),
and substituting this ansatz into (\ref{quasi}), we obtain 
\begin{equation}
\partial_ib=\mu^{i}\partial_ia, ~~ \partial_ip=\lambda^{i}\partial_ia, ~~ \partial_iq=\mu^{i}\lambda^i\partial_ia, ~~ \partial_ic=(\mu^{i})^2\partial_ia,
\label{bpqc}
\end{equation}
$\partial_i=\partial/\partial_{R^i}$, along with the dispersion relation
\begin{equation}
D(\lambda_i, \mu_i) =f_{a}+f_{b}\mu^{i}+f_{c}(\mu^{i})^2+f_{p}\lambda^{i}+f_{q}\lambda^{i}\mu^{i} - (\lambda^{i})^2=0;
\label{disp}
\end{equation}
in what follows we assume that the dispersion relation (\ref{disp}) defines an irreducible conic in the $(\lambda, \mu)$-plane. This is equivalent to the requirement that the expression 
$$
\Delta = f_b^2+f_bf_pf_q-f_af_q^2-f_cf_p^2-4f_af_c
$$ 
does not vanish. 
 Calculating the consistency conditions for (\ref{bpqc}) we obtain  
\begin{equation}
\begin{array}{c}
\partial _{i}\partial _{j}a=\frac{\partial _{j}\lambda ^{i}}{\lambda
^{j}-\lambda ^{i}}\partial _{i}a+\frac{\partial _{i}\lambda ^{j}}{\lambda
^{i}-\lambda ^{j}}\partial _{j}a, \\
\ \\
\partial_j\lambda^i\partial_ia+ \partial_i\lambda^j\partial_ja=0.
\end{array}
 \label{Comp}
\end{equation}
Applying the operator $\partial_j$ to the dispersion relation (\ref{disp}) and using (\ref{bpqc}) and (\ref{Com}), we obtain  $\partial_j \lambda^i$ and $\partial_j \mu^i$ in the form
\begin{equation}
\partial_j \lambda^i=(\lambda^i-\lambda^j)B_{ij}\partial_ja, \ \ \ \ \
\partial_j \mu^i=(\mu^i-\mu^j)B_{ij}\partial_ja
\label{lm}
\end{equation}
where $B_{ij}$ are rational expressions in $\lambda^i, \lambda^j, \mu^i, \mu^j$ whose coefficients depend on the second order partial derivatives of the function $f(a,b,c,p,q)$. Explicitly, one has
$$
B_{ij}=\frac{N_{ij}}{D_{ij}}
$$
where
\begin{equation}
\begin{array}{c}
N_{ij}=f_{aa}+f_{ab}(\mu^i+\mu^j)+f_{ac}\big((\mu^i)^2+(\mu^j)^2\big)+f_{ap}(\lambda^i+\lambda^j)+f_{aq}(\lambda^i\mu^i+\lambda^j\mu^j)\\
\ \\
+f_{bb}\mu^i\mu^j+f_{bc}\mu^i\mu^j(\mu^i+\mu^j)+f_{bp}(\lambda^i\mu^j+\lambda^j\mu^i)+f_{bq}\mu^i\mu^j(\lambda^i+\lambda^j) \\
\ \\
+f_{cc}(\mu^i)^2(\mu^j)^2+f_{cp}\big(\lambda^i(\mu^j)^2+\lambda^j(\mu^i)^2\big)+f_{cq}\mu^i\mu^j(\lambda^i\mu^j+\lambda^j\mu^i) \\
\ \\
+f_{pp}\lambda^i\lambda^j+f_{pq}\lambda^i\lambda^j(\mu^i+\mu^j)+f_{qq}\lambda^i\lambda^j\mu^i\mu^j
\end{array}
\end{equation}
and 
\begin{eqnarray}
D_{ij}&=&
-2\lambda^i\lambda^j+2f_a+f_b(\mu^i+\mu^j)+2f_c\mu^i\mu^j+f_p(\lambda^i+\lambda^j)+f_q(\lambda^i\mu^j+\lambda^j\mu^i) \nonumber \\
&=&
4D\left(\frac{\lambda^i+\lambda^j}{2}, \ \frac{\mu^i+\mu^j}{2}\right),
\end{eqnarray}
the last equality holding modulo the dispersion relation (\ref{disp}).
Substituting (\ref{lm}) into $(\ref{Comp})_1$ one obtains
\begin{equation}
\partial_i\partial_j a=-2B_{ij}\partial_ia \partial_j a;
\label{aij}
\end{equation}
 the second condition $(\ref{Comp})_2$ holds identically by virtue of the symmetry $B_{ij}=B_{ji}$. 
The compatibility conditions $\partial_k\partial_j
\lambda^i=\partial_j\partial_k \lambda^i$,
 $\partial_k\partial_j \mu^i=\partial_j\partial_k \mu^i$ and
$\partial_k\partial_j\partial_i a=\partial_j\partial_k\partial_i a$ are
equivalent to the relations
\begin{equation}
\partial_kB_{ij}=(B_{ij}B_{kj}+B_{ij}B_{ik}-B_{kj}B_{ik})\partial_ka,
\label{int}
\end{equation}
which must be satisfied identically by virtue of (\ref{bpqc}), (\ref{disp}) and
(\ref{lm}).

{\bf Remark.} Notice that each compatibility condition involves three distinct  indices only. This observation immediately implies that

\noindent {\bf (i)} any equation of the form (\ref{main}) possesses infinitely many two-component reductions parametrized by two arbitrary functions of one variable. Indeed, these reductions are governed by the equations (\ref{lm}) and (\ref{aij}) where $i, j=1, 2$ and $\lambda^i, \mu^i$ satisfy the dispersion relation (\ref{disp}). These equations are automatically consistent, and the general solution depends, modulo reparametrizations $R^i\to \varphi^i(R^i)$, on two arbitrary functions of one variable.
Therefore, the existence of two-component reductions is a common phenomenon which is not related to the integrability.

\noindent {\bf (ii)}  the existence of three-component reductions implies the existence of $n$-component reductions for arbitrary $n$. Thus, one can define the integrability as the existence of infinitely many three-component reductions parametrized by three arbitrary functions of one variable. This property is reminiscent of the well-known three-soliton condition in the theory of integrable systems. 

In order to simplify the derivation of the integrability conditions we first  rewrite (\ref{int}) as
\begin{equation}
\label{Nij}
\partial_kN_{ij}=N_{ij}\left(\frac{1}{D_{ij}}\partial_kD_{ij} +
B_{kj}\partial_ka+B_{ik}\partial_ka\right)-D_{ij}B_{kj}B_{ik}\partial_ka.
\end{equation}
Here, the third order derivatives of $f(a,b,c,p,q)$ are present only in
the l.h.s.\ term $\partial_kN_{ij}$ (see the form of $N_{ij}$ and $ D_{ij}$ above). Further reduction of the complexity
of the expression in the r.h.s.\ is achieved by representing
$1/D_{ij}$ in the form
$$
\begin{array}{c}
\displaystyle \frac{1}{D_{ij}} = U_{ij}=
[2 f_a  + (\mu_1+\mu_2) f_b + 2 \mu_1 \mu_2 f_c - (\lambda_1 + \lambda_2) f_p - 
(\lambda_1 \mu_2 + \lambda_2 \mu_1) f_q \\
\ \\
+ f_p^2+ (\mu_1 + \mu_2) f_p f_q + \mu_1 \mu_2 f_q^2 + 2 \lambda_1 \lambda_2]/((\mu_1-\mu_2)^2 \Delta)
\end{array}
$$
(which holds identically modulo the dispersion relation (\ref{disp})),
and the subsequent substitution $B_{st}=N_{st}/D_{st} = N_{st}U_{st}$.
The denominators of the r.h.s.\ terms in (\ref{Nij}) cancel out,  producing
a polynomial in  $\lambda^i$, $\lambda^j$, $\lambda^k$, $\mu^i$,
$\mu^j$, $\mu^k$ with coefficients depending on the derivatives of the
density $f(a,b,c,p,q)$.
This was the crucial simplification of the calculation:
the starting expression for the r.h.s.\ of (\ref{Nij}) has more than 450.000 
terms
with different denominators;
after  properly organized cancellations we get a polynomial expression with less
than 6.000 terms.
Using the dispersion relation (\ref{disp}), we simplify this polynomial by substituting
the powers of $(\lambda^i)^s$,
$(\lambda^j)^s$, $(\lambda^k)^s$,  $s\geq 2$,
arriving at a polynomial of degree one in each of $\lambda^i$,
$\lambda^j$, $\lambda^k$, and degree two in $\mu$'s.
Equating similar coefficients  in both sides of (\ref{Nij}),
we arrive at a set of 35 equations for the derivatives of the function
$f(a, b, c, p, q)$, which are linear in the third order derivatives.
Solving this linear system we obtain the  closed form expressions for all third order derivatives
of $f(a,b,c,p,q)$ in terms of its first and second order derivatives, which we represent  symbolically in the form
\begin{equation}
d^3f=R(df, d^2f);
\label{*}
\end{equation}
here $R$ denotes complicated rational expressions in the first and second order derivatives of $f$ (it is clear from the above that the function $f$ itself will not enter these expressions explicitly). 
 A straigtforward computer calculation shows that the overdetermined system (\ref{*}) is in involution. Thus, the moduli space of integrable systems of the type (\ref{main}) is $21$-dimensional: one can arbitrarily specify the values of $f$, $df$ and $d^2f$ at any fixed point. This amounts to $1+5+15=21$ arbitrary constants. This finishes the proof of Theorem 1.
 
 The right hand sides of (\ref{*})   are not presented here because of their complexity (see Sect. 3 for a discussion of particular cases where these formulae become less cumbersome, and are presented explicitly). The integrability conditions  (\ref{*}) provide a straightforward computer test of integrability for any  equation from the class under consideration, and allow one to obtain classification results. These conditions also allow computer checks of various  properties of differential-geometric objects  naturally associated with  equations from the class (\ref{main}). 
 

{\bf Remark.} In two dimensions, any second order PDE of the form $F(u_{xx}, u_{xy}, u_{yy})=0$ is automatically integrable. Indeed, introducing the parametrization $u_{xx}=a, \ u_{xy}=b, \ u_{yy}=f(a, b)$, so that $F(a, b, f(a, b))=0$, one obtains a two-component quasilinear system 
$$
a_y=b_x, ~~~ b_y=f(a, b)_x;
$$
any system of this type linearises under a hodograph transformation which interchanges dependent and independent variables. This simple trick, however,  does not work in more than two dimensions.

\section{Examples and classification results}

In this section we discuss  examples (both known and new) and partial classification results of integrable equations of the form (\ref{main}) which are obtained using the integrability conditions (\ref{*}) derived in Sect. 2. 

\subsection{Integrable equations of the form $u_{tt}=f(u_{xx}, u_{yy})$}

\medskip

\noindent
In this case the  integrability conditions (\ref{*}) simplify to
\begin{eqnarray}
&&f_{aaa}=f_{aa}\left(\frac{f_{ac}}{f_c}+\frac{f_{aa}}{f_a}\right),\quad
f_{aac}=f_{aa}\left(\frac{f_{cc}}{f_c}+\frac{f_{ac}}{f_a}\right),\notag \\
&&f_{acc}=f_{cc}\left(\frac{f_{cc}}{f_c}+\frac{f_{ac}}{f_a}\right),\quad
f_{ccc}=f_{cc}\left(\frac{f_{cc}}{f_c}+\frac{f_{ac}}{f_a}\right),\quad
f_{aa}f_{cc}=(f_{ac})^2;\notag
\end{eqnarray}
recall that, according to the notation introduced in Sect. 2, $a=u_{xx}, ~ c=u_{yy}$. This system is in involution and its general solution depends on 5 integration constants. To solve it explicitly one notices that the first two equations imply  $\frac{f_{aa}}{f_af_c}=const$. Similarly, the next two equations imply  $\frac{f_{cc}}{f_af_c}=const$. Further elementary integration leads to the following general solution,
 $$
 f(a, c)=s\ln(me^{\mu a}+ ne^{\nu c})+const.
 $$
 This corresponds to  equations of the form
 $$
 me^{\mu u_{xx}}+ ne^{\nu u_{yy}}+ ke^{\kappa u_{tt}}=0,
 $$
 as well as  degenerations thereof,
  $$
 mu_{xx}+ nu_{yy}+ ke^{\kappa u_{tt}}=0,
 $$
where the coefficients are arbitrary constants. In fact, all these coefficients can be eliminated by appropriate complex 
rescalings, leading to the  two essentially different examples,
$$
e^{u_{xx}}+ e^{u_{yy}}= e^{u_{tt}} \ \ \ \  {\rm and} \ \ \ \  u_{xx}+ u_{yy}= e^{u_{tt}}.
$$
The first equation is apparently new, while the second is the Boyer-Finley equation.

\subsection{Integrable equations of the form $u_{tt}=f(u_{xx}, u_{xy}, u_{yy})$}

\medskip

This  is a generalization  of the previous case. Thus, we assume $f_b\ne 0$. The integrability conditions (\ref{*}) take the form
\begin{eqnarray*}
&&f_{aa}f_{bb}={f_{ab}}^2,\quad
f_{aa}f_{cc}={f_{ac}}^2, \quad
f_{bb}f_{cc}={f_{bc}}^2,  \notag \\
\\
&&f_{aa}f_{bc}=f_{ab}f_{ac}, \quad
f_{ab}f_{cc}=f_{ac}f_{bc}, \quad 
f_{ab}f_{bc}=f_{ac}f_{bb},  \notag \\
\\
&&f_{aaa}=\frac{2f_{aa}\left(f_bf_{ab}-2f_cf_{aa}-2f_af_{ac}\right)}{{f_b}^2-4f_af_c}, \quad
f_{aab}=\frac{2f_{ab}\left(f_bf_{ab}-2f_cf_{aa}-2f_af_{ac}\right)}{{f_b}^2-4f_af_c},\notag \\
\\ 
&&f_{aac}=\frac{2f_{ac}\left(f_bf_{ab}-2f_cf_{aa}-2f_af_{ac}\right)}{{f_b}^2-4f_af_c}, \quad
f_{abb}=\frac{2f_{ab}\left(f_bf_{bb}-2f_cf_{ab}-2f_af_{bc}\right)}{{f_b}^2-4f_af_c}, \notag \\
\\
&&f_{acc}=\frac{2f_{ac}\left(f_bf_{bc}-2f_cf_{ac}-2f_af_{cc}\right)}{{f_b}^2-4f_af_c},  \quad
f_{abc}=\frac{2f_{ab}\left(f_bf_{bc}-2f_cf_{ac}-2f_af_{cc}\right)}{{f_b}^2-4f_af_c},\notag \\
\\
&&f_{bbb}=\frac{2f_{bb}\left(f_bf_{bb}-2f_cf_{ab}-2f_af_{bc}\right)}{{f_b}^2-4f_af_c},\quad
f_{bbc}=\frac{2f_{bc}\left(f_bf_{bb}-2f_cf_{ab}-2f_af_{bc}\right)}{{f_b}^2-4f_af_c},\notag \\
\\
&&f_{bcc}=\frac{2f_{bc}\left(f_bf_{bc}-2f_cf_{ac}-2f_af_{cc}\right)}{{f_b}^2-4f_af_c},\quad
f_{ccc}=\frac{2f_{cc}\left(f_bf_{bc}-2f_cf_{ac}-2f_af_{cc}\right)}{{f_b}^2-4f_af_c}.  \notag \\
\end{eqnarray*}
Notice that the first six equations imply that the Hessian matrix of the function $f(a, b, c)$ has rank one,
so that one can set $f_a=m(f_b), \ f_c=n(f_b)$. The substitution into the remaining equations implies either $m''=n''=0$,  or $f_{bb}=0$. Further elementary analysis leads, up to linear transformations of $x$ and $y$, to the following canonical forms (here we only list those representatives which contain a nontrivial dependence on $u_{xy}$):

\noindent {\bf (a)} $u_{tt}=\alpha u_{xx}+\beta u_{yy}+\varphi (u_{xy})$, here $\varphi$ satisfies a third order ODE $\varphi '''((\varphi')^2-4\alpha \beta)=2\varphi'(\varphi'')^2$. The integration  leads to the three canonical forms,
$$
\begin{array}{c}
u_{tt}=\alpha u_{xx}+\beta u_{yy}+\frac{2\sqrt {\alpha \beta}}{\gamma}\ln \cosh \gamma u_{xy}, \\
\ \\
u_{tt}=\alpha u_{xx}+\gamma \ln u_{xy}, ~~~~~ {\rm or} ~~~~~ u_{tt}= \ln u_{xy}, \\
\end{array}
$$
the last two cases corresponding to $\beta=0$ and $\alpha=\beta=0$, respectively. The first two examples can be viewed as generalizations of the Boyer-Finley equation.

\noindent {\bf (b)} $u_{tt}=u_{xy}+\beta (u_{xx}-u_{yy})+\varphi (u_{xx}+u_{yy})$, here $\varphi '''(4(\varphi')^2-1-4 \beta^2)=8\varphi'(\varphi'')^2$. This leads to the equation
$$
u_{tt}=u_{xy}+\beta(u_{xx}-u_{yy})+\frac{\sqrt{1+4\beta^2}}{2\gamma}\ln \cosh \gamma(u_{xx}+u_{yy}).
$$

\noindent {\bf (c)}  $u_{tt}= u_{xy}+\beta u_{xx}+\varphi (u_{yy})$, here $\varphi$ satisfies a third order ODE $\varphi '''(4\beta \varphi'-1)=4\beta(\varphi'')^2$. This results in the equation
$$
u_{tt}=u_{xy}+\beta u_{xx}+\frac{1}{4\beta}u_{yy}+\alpha e^{\gamma u_{yy}},
$$
whose degeneration, corresponding to $\beta=0$, is the dKP equation $u_{tt}=u_{xy}+u_{yy}^2$.

\subsection{Integrable equations of the form $u_{xy}=f(u_{xt}, u_{yt})$}

\medskip

Formally,  equations from this class are not of  the  form (\ref{evol}), however, they can readily be made `evolutionary'  by an appropriate linear change of the independent variables. The resulting set of integrability conditions looks as follows:
\begin{eqnarray}
&&f_{ppp}=f_{pp}\left(\frac{f_{pq}}{f_q}+\frac{f_{pp}}{f_p}\right),\quad
f_{ppq}=f_{pp}\left(\frac{f_{qq}}{f_q}+\frac{f_{pq}}{f_p}\right),\notag \\
&&f_{pqq}=f_{qq}\left(\frac{f_{pq}}{f_q}+\frac{f_{pp}}{f_p}\right),\quad
f_{qqq}=f_{qq}\left(\frac{f_{qq}}{f_q}+\frac{f_{pq}}{f_p}\right).\notag
\end{eqnarray}
This system is in involution, and its general solution depends on 6 integration constants. To solve it explicitly one notices that the first two equations imply  $\frac{f_{pp}}{f_pf_q}=const$. Similarly, the last two equations imply  $\frac{f_{qq}}{f_pf_q}=const$. Further elementary integration gives, under the assumption that both $f_{pp}$ and $f_{qq}$ are nonzero, the  general solution
 $$
f(p, q)=\frac{1}{\kappa}\ln\left(\frac{me^{\mu p}+ ne^{\nu q}}{re^{\mu p+\nu q}-k}\right),
$$
which leads to the dipersionless Hirota-type equation for the BKP hierarchy \cite{Bogdanov},
 \begin{equation}
 me^{\mu u_{xt}}+ ne^{\nu u_{yt}}+ ke^{\kappa u_{xy}}=re^{\mu u_{xt}+\nu u_{yt}+\kappa u_{xy}}.
 \label{Hirota}
\end{equation}
Up to complex rescalings, it can be transformed to 
 $$
 e^{ u_{xt}}+ e^{ u_{yt}}+ e^{ u_{xy}}=e^{ u_{xt}+ u_{yt}+ u_{xy}}.
 $$
 Degenerations of (\ref{Hirota}), corresponding to  $f_{pp}=0$ or $f_{pp}=f_{qq}=0$, result in
 $$
u_{xy}= u_{xt}+e^{u_{yt}}, \ \ \ \ \ \  u_{xy}= u_{xt}\tan(u_{yt})
 $$
and
$$
u_{xy}=u_{xt}u_{yt},
$$
respectively.

\subsection{Integrable equations of the form $u_{tt}=f(u_{xx},u _{xt},u _{xy})$}

\medskip

Equations of this type arise in the theory of  integrable hydrodynamic chains satisfying the additional`Egorov' property. 
Here we reproduce  the classification result from \cite{MaksEgor} (see also \cite{Buch}, \cite{Fer4}). The integrability conditions (\ref{*}) take the form
\begin{eqnarray}
&&f_{bbb}=\frac{2f_{bb}^{2}}{f_{b}},\quad
f_{abb}=\frac{2f_{ab}f_{bb}}{f_{b}}, \quad
f_{pbb}=\frac{2f_{pb}f_{bb}}{f_{b}},  \notag \\
&&f_{aab}=\frac{2f_{ab}^{2}}{f_{b}}, \quad
f_{apb}=\frac{2f_{ab}f_{pb}}{f_{b}}, \quad
f_{ppb}=\frac{2f_{pb}^{2}}{f_{b}},  \notag \\
&&f_{ppp}=\frac{2}{f_{b}^{2}}\left(
f_{p}f_{pb}^{2}+f_{pb}(f_{b}f_{pp}+2f_{ab})-f_{bb}(f_{p}f_{pp}+2f_{ap})
\right),  \notag \\
&&f_{app}=\frac{2}{f_{b}^{2}}\left(
f_{a}f_{pb}^{2}+f_{ab}(f_{b}f_{pp}+f_{ab})-f_{bb}(f_{a}f_{pp}+f_{aa})\right), \\
&&f_{aap}=\frac{2}{f_{b}^{2}}\left(
f_{bb}(f_{p}f_{aa}-2f_{a}f_{ap})-f_{ab}(f_{p}f_{ab}-2f_{b}f_{ap})-f_{pb}(f_{b}f_{aa}-2f_{a}f_{ab})\right),
\notag \\
&&f_{aaa}=\frac{2}{f_{b}^{2}}\left(
(f_{a}+f_{p}^{2})f_{ab}^{2}+f_{a}^{2}f_{pb}^{2}+f_{b}^{2}(f_{ap}^{2}-f_{aa}f_{pp})-f_{pp}f_{bb}f_a^2 \right.
\notag \\
&&\qquad
+f_{ab}f_{b}(f_{aa}+2(f_{a}f_{pp}-f_{p}f_{ap}))+2f_{pb}(f_{p}(f_{b}f_{aa}-f_{a}f_{ab})-f_{a}f_{b}f_{ap})
\notag \\
&&\qquad -\left.
f_{bb}((f_{a}+f_{p}^{2})f_{aa}-2f_{a}f_{p}f_{ap})\right); \notag
\end{eqnarray}
this system is in involution and its general solution depends on $10 $
arbitrary constants.

The integration of these equations leads to the four essentially different canonical forms,
\begin{eqnarray*}
u_{tt} &=&u_{xy}+\frac{1}{4A}(Au_{xt}+2Bu_{xx})^{2}+Ce^{-Au_{xx}}, \\
u_{tt} &=&\frac{u_{xy}}{
u_{xx}}+\left( \frac{1}{u_{xx}}
+\frac{A}{4u_{xx}^{2}}\right) u_{xt}^{2}+\frac{B}{u_{xx}^{2}}u_{xt}+\frac{B^{2}}{Au_{xx}^{2}}+Ce^{A/u_{xx}},
\\
u_{tt} &=&\frac{u_{xy}}{u_{xt}}+\frac{1}{6}\eta (u_{xx})u_{xt}^{2}, \\
u_{tt} &=&\ln u_{xy}-\ln \theta _{1}\left(u_{xt}, u_{xx}\right) -\frac{1}{4}\overset{u_{xx}}{\int }\eta (\tau )d\tau,
\end{eqnarray*}
see \cite{MaksEgor}. Here $A, B, C$ are arbitrary constants,  $\eta $ is a solution to the Chazy equation  \cite{Chazy},
$$
\eta^{\prime \prime \prime }+2\eta \eta ^{\prime \prime }=3(\eta^{\prime}) ^2,
$$
and $\theta_1$ is a theta-function in the variable $u_{xt}$ whose dependence on $u_{xx}$ is governed by the Chazy equation.





\subsection{Symplectic Monge-Amp\'ere equations}

Let us consider a function $u(x^1, ... x^k)$ of $k$ independent variables, and introduce the $k\times k$ Hessian matrix   $U=[ u_{ij} ]$ of its second order partial derivatives. The symplectic Monge-Amp\'ere equation is a PDE of the form
$$
M_k+M_{k-1}+... +M_1+M_0=0
$$
where $M_l$ is a constant-coefficient linear combination of all  $l\times l$ minors of  the matrix $U$, $0\leq l\leq k$. Thus,  $M_k= det~U= Hess~u$,  $M_0$ is a constant, etc. Equivalently, these PDEs can be obtained by equating to zero a constant-coefficient  $k$-form in the $2k$ variables $x^i, u_i$. This class of equations is invariant under the natural action  of the symplectic group $Sp(2k)$. In the case $k=2$ one obtains a standard Monge-Amp\'ere equation,
\begin{equation}
u_{11}u_{22}-u^2_{12}+\alpha u_{11}+\beta u_{12}+\gamma u_{22}+\delta =0,
\label{MAA}
\end{equation}
which can be interpreted as the equation of a `sphere' corresponding to the pseudo-Euclidean metric $du_{11}du_{22}-du^2_{12}$.
Monge-Amp\'ere equations (\ref{MAA}) can be characterized as the only equations of the form $F(u_{11}, u_{12}, u_{22})=0$ which are linearizable by a transformation from the equivalence group $Sp(4)$. 

The case $k=3$ is also understood completely: one can show that, for $k=3$, any symplectic Monge-Amp\'ere equation is either linearizable  (in this case it is automatically integrable), or $Sp(6)$-equivalent to either of the  two essentially different canonical forms \cite{Lychagin, Banos},
\begin{equation}
Hess~u=1, ~~~~~  Hess~u=u_{11}+u_{22}+u_{33}.
\label{MA}
\end{equation}
Based on the integrability conditions (\ref{*}), we have verified directly  that  both  PDEs (\ref{MA}) are {\it not} integrable by the method of hydrodynamic reductions. Thus, a 3-dimensional symplectic Monge-Amp\'ere equation is integrable if and only if it is linearizable (this is no longer true in more than three dimensions - see Sect. 5). The linearizability condition constitutes a single relation among the coefficients of the equation: for a Monge-Amp\'ere equation of the form
\begin{equation}
\begin{array}{c}
\epsilon \ det \left[
\begin{array}{ccc}
u_{11} & u_{12} & u_{13} \\
u_{12} & u_{22} & u_{23} \\
u_{13} & u_{23} & u_{33} 
\end{array}
\right]
+h_1(u_{22}u_{33}-u_{23}^2)+h_2(u_{11}u_{33}-u_{13}^2)+h_3(u_{11}u_{22}-u_{12}^2) \\
\ \\
+g_1(u_{11}u_{23}-u_{12}u_{13})+g_2(u_{22}u_{13}-u_{12}u_{23})+g_3(u_{33}u_{12}-u_{13}u_{23}) \\
\ \\
+s_1 u_{11}+s_2 u_{22}+s_3 u_{33}+\tau_1 u_{23}+\tau_2 u_{13} +\tau_3 u_{12}+\nu=0,
\end{array}
\label{EQ}
\end{equation}
the linearizability condition defines a  quartic hypersurface in the space of coefficients,
\begin{equation}
\begin{array}{c}
h_1^2s_1^2+h_2^2s_2^2+h_3^2s_3^2+g_1^2s_2s_3+g_2^2s_1s_3+g_3^2s_1s_2 \\
\ \\
-2(h_1h_2s_1s_2+h_1h_3s_1s_3+h_2h_3s_2s_3)+4\epsilon s_1s_2s_3+4\nu h_1h_2h_3 \\
\ \\
+\epsilon \tau_1 \tau_2 \tau_3-\nu g_1g_2g_3-\epsilon^2 \nu^2-\nu (g_1^2h_1+g_2^2h_2+g_3^2h_3)\\
\ \\
-(g_1\tau_1+g_2\tau_2+g_3\tau_3+2\epsilon \nu)(h_1s_1+h_2s_2+h_3s_3-\epsilon \nu) \\
\ \\
+2(g_1h_1s_1\tau_1+g_2h_2s_2\tau_2+g_3h_3s_3\tau_3)+\tau_1^2h_2h_3+\tau_2^2h_1h_3+\tau_3^2h_1h_2\\
\ \\
-\epsilon(\tau_1^2s_1+\tau_2^2s_2+\tau_3^2s_3)+s_1\tau_1g_2g_3+s_2\tau_2g_1g_3+s_3\tau_3g_1g_2\\
\ \\
-(g_1h_1\tau_2\tau_3+g_2h_2\tau_1\tau_3+g_3h_3\tau_1\tau_2)=0.
\end{array}
\label{quartic}
\end{equation}
Notice that, if $\epsilon \ne 0$,  one can always eliminate second order minors in (\ref{EQ}) by adding to $u$ an appropriate quadratic form. In this case the equation takes the form
$$
det \left[
\begin{array}{ccc}
u_{11} & u_{12} & u_{13} \\
u_{21} & u_{22} & u_{23} \\
u_{31} & u_{32} & u_{33} 
\end{array}
\right]
+s_1 u_{11}+s_2 u_{22}+s_3 u_{33}+\tau_1 u_{23}+\tau_2 u_{13} +\tau_3 u_{12}+\nu=0,
$$
while  the linearizability condition simplifies to
$$
4s_1s_2s_3+\nu^2+\tau_1 \tau_2 \tau_3 -s_1 \tau_1^2-s_2 \tau_2^2-s_3 \tau_3^2=0.
$$
\noindent {\bf Remark 1.} The condition (\ref{quartic}) can be given an invariant formulation \cite{Lychagin, Banos}. Let $\Omega=du_i\wedge dx^i$ be the standard symplectic form, and $v=v_{\Omega}$ the associated canonical bivector. Let $\omega$ be the effective 3-form corresponding to the equation (\ref{EQ}): recall that $\omega$ is effective if $\omega \wedge \Omega =0$. Let us  define the quadratic form
\begin{equation}
\bot^2(i_X\omega \wedge i_X\omega)
\label{quadratic}
\end{equation}
where $i_X\omega$ is the inner product of $\omega$ with a vector $X$, and $\bot$ is the operator of inner multiplication by the canonical bivector. The linearizability condition (\ref{quartic}) is equivalent to the requirement that the quadratic form (\ref{quadratic}) is degenerate.

\noindent {\bf Remark 2.}  The equation $(\ref{MA})_1$ arises in the theory of improper affine spheres, while the equation $(\ref{MA})_2$ describes special Lagrangian $3$-folds in $C^3$ \cite{Calabi, Joyce}. Their contact-equivalent forms  are
$$
u_{tt}=u_{xx}u_{yy}-u_{xy}^2
$$ 
and
$$
u_{xx}u_{yy}-u_{xy}^2+u_{xx}u_{tt}-u_{xt}^2+u_{tt}u_{yy}-u_{ty}^2=1,
$$
respectively (both non-integrable).

\noindent {\bf Remark 3.} Equations of the form (\ref{EQ}) and the linearizability condition (\ref{quartic}) have a clear projective-geometric interpretation. 
Let us consider the Lagrangian Grassmannian $\Lambda$, which can be  (locally) identified  with the space of $3\times 3$ symmetric matrices $U$ (see Sect. 4 for more details). The minors of $U$ define the Veronese embedding of $\Lambda$ into the projective space $P^{13}$ (we will identify $\Lambda$ with its projective embedding). Thus,  symplectic Monge-Amp\'ere equations  correspond to  hyperplane sections of $\Lambda$. Linearizable equations correspond to hyperplanes which are tangential to  $\Lambda$. Thus, the linearizability condition (\ref{quartic}) coincides with the equation of the dual variety to $\Lambda$, which is known to be a hypersurface of degree four  in $(P^{13})^*$ (see  \cite{Landsberg1, Landsberg2, Mukai} for a general theory behind this example).

\medskip

The results of this section answer in the negative the question  formulated by Joyce in \cite{Joyce}:
{\it do special Lagrangian $m$-folds in $C^m$ for $m\geq 3$ constitute  some kind of higher-dimensional integrable system}?

\section{Differential-geometric aspects of integrable equations of the dispersionless Hirota type}

In this section we adopt a geometric point of view, and consider  (\ref{main}) as the equation difining a $5$-dimensional hypersurface $M^5$ in the $6$-dimensional space with coordinates $u_{ij}$ (which, as we will see shortly, is naturally identified with the Lagrangian Grassmannian $\Lambda$). Our aim  is to reformulate the integrability conditions in the intrinsic geometric terms. 

\medskip

\noindent In Sect. 4.1 we discuss differential geometry of the Lagrangian Grassmannian $\Lambda$. We introduce the main object needed to develop differential geometry of hypersurfaces, namely,   a canonically defined symmetric cubic form  which is invariant under the natural action of the symplectic group $Sp(6)$. This form (to be precise, its conformal class), plays the role of an ambient  flat conformal structure. 

\medskip

\noindent In Sect. 4.2 we develop differential geometry of hypersurfaces  $M^5\subset \Lambda$. The induced cubic form gives rise to  a field of rational normal curves of degree four defined in the projectivization of the tangent bundle  $TM^5$. Thus, a hypersurface of the Lagragian Grassmannian carries an intrinsic  irreducible $GL(2)$-structure.

\medskip

\noindent In Sect. 4.3 we show that two-component and three-component hydrodynamic reductions of the equation (\ref{main}) correspond to bisecant  and trisecant submanifolds of $M^5$, respectively. This allows one to reformulate  the integrability  as the existence of sufficiently many trisecant submanifolds parametrized by three arbitrary functions of one variable.

\medskip

\noindent Some useful computational formulas are provided in Sect. 4.4.

\medskip

\noindent General aspects of  irreducible $GL(2)$-structures in dimension $5$ are discussed in Sect. 4.5. We demonstrate that the requirement of the existence of  bisecant surfaces imposes strong constraints on the torsion of the corresponding $GL(2)$-structure.

\subsection{$Sp(6)$ as a symmetry group of the problem: differential geometry of the Lagrangian Grassmannian}

Let us consider a $6$-dimensional symplectic space with canonical coordinates ${\bf x}=(x^1, x^2, x^3)^t$ and ${\bf p}=(p_1, p_2, p_3)^t$, viewed as $3$-component column vectors.  Lagrangian submanifolds can be  parametrized in terms of a generating function $u(x, y, t)$: ${\bf x}=(x, y, t)^t$ and
${\bf p}=(u_x, u_y, u_t)^t$.
Lagrangian planes are defined by the equation $d{\bf p}=U d{\bf x}$ where $U$ is a $3\times 3$ symmetric matrix (the Hessian matrix of $u$). Thus, the Lagrangian Grassmannian $\Lambda$ is $6$-dimensional, and can (locally)  be identified with the  space of $3\times 3$ symmetric matrices. The equation (\ref{main}) defines a $5$-dimensional hypersurface $M^5\subset \Lambda$; the corresponding solutions $u(x, y, t)$ can be interpreted as Lagrangian submanifolds whose Gaussian images belong to the hypersurface $M^5$. 

The  action of the linear symplectic group $Sp(6)$,
$$
\left(\begin{array}{c}
d{\bf \tilde p}
 \\
 d{\bf \tilde x}
 \end{array}\right)=
\left( \begin{array}{cc}
 A & B
 \\
C & D 
 \end{array}\right)
\left( \begin{array}{c}
d{\bf  p}
 \\
 d{\bf  x}
 \end{array}\right),
 $$
 naturally extends to $\Lambda$:
 \begin{equation}
 \tilde U=(AU+B)(CU+D)^{-1}.
 \label{Sp}
 \end{equation}
Here $A, B, C, D$ are $3\times 3$ matrices such that $A^tC=C^tA, \  B^tD=D^tB, \ A^tD-C^tB=id$; notice that  the extended action is no longer linear. The transformation law (\ref{Sp}) suggests that the action of $Sp(6)$ preserves the class of equations (\ref{main}), indeed, second order derivatives transform through  second order derivatives only. Moreover, since any changes of variables obviously preserve the integrability, the group $Sp(6)$ can be viewed as a natural {\it equivalence group} of the problem: it maps integrable equations to integrable. Thus,  $Sp(6)$ is a point symmetry group of the integrability conditions (\ref{*}) derived in Sect. 2. The classification of integrable equations of the form (\ref{main}) has to be performed modulo this equivalence: two $Sp(6)$-related equations should be regarded as `the same'. 

Geometrically, our problem is reduced to the classification of hypersurfaces $M^5$ of the Lagrangian Grassmannian $\Lambda$ (satisfying certain `integrability' conditions  to be specified later, see Sect. 4.3), up to the action of $Sp(6)$. Our next goal is to clarify which differential-geometric structures are  induced on  hypersurfaces of $\Lambda$. To do so one first needs to introduce differential-geometric objects which are naturally defined on the Lagrangian Grassmannian, and are invariant under the action of $Sp(6)$.

The condition of nontrivial intersection of  two infinitesimally close Lagrangian planes  defined by symmetric matrices $U$ and $U+dU$, is  $\det dU=0$. This condition is manifestly invariant under the action of $Sp(6)$ as specified by (\ref{Sp}). Thus, each tangent space to the Lagrangian Grassmannian $\Lambda$ is equipped with a cubic cone $C$ defined by the equation $\det dU=0$. The projectivisation of this  cone is known as the `cubic symmetroid'. The singular locus $V\subset C$ of the cubic cone  is specified by the condition $rank~ dU=1$. Its projectivization is known as the Veronese variety, which  is a non-singular algebraic surface of degree four.

The identity
$$
\det d\tilde U=\frac{\det (A-(AU+B)(CU+D)^{-1}C) }{\det (CU+D)}\ \det dU,
$$
which readily follows from (\ref{Sp}), implies that the conformal class of the cubic form $\det dU$ is invariant under the action of $Sp(6)$. The converse is also true:

\begin{theorem} The group of conformal automorphisms of the symmetric cubic form $\det dU$ is isomorphic to $Sp(6)$.
\end{theorem} 
 The proof consists of a direct calculation of conformal automorphisms of the cubic form  $\det dU$. The corresponding infinitesimal generators are presented in Sect. 6. In fact,  Theorem 2 can be traced back to the old paper by Cartan \cite{Cartan}, Theorem XX. We  refer to \cite{Bertram, Gindikin1} for further generalizations of this Liouville-type result. 

{\bf Remark.} The  theorem is true in any dimension. Thus, for the Grassmannian of  Lagrangian 2-planes in  a 4-dimensional symplectic space (naturally identified with  $2\times 2$ symmetric matrices $U$), the expression $\det dU$ defines  the Lorentzian metric $du_{11}du_{22}-du_{12}^2$ whose  conformal automorphisms generate the group $SO(3, 2)$. This establishes a well-known isomorphism between $Sp(4)$ and $SO(3, 2)$.

\subsection{Geometric structures on a hypersurface of the Lagrangian Grassmannian}


Let $M^5$ be a hypersurface of the Lagrangian Grassmannian $\Lambda$. Taking a point $s\in M^5$ and intersecting the tangent space $T_sM^5$ with the singular locus $V$ of the cubic cone $C$ in $T_s\Lambda$, one obtains, after a projectivisation,   a rational normal curve of degree four (it is  well-known  that a generic hyperplane section of the Veronese variety $V$ is a rational normal curve). This curve can also be interpreted as a set of matrices of rank one in the tangent space $T_sM^5$ (recall that $T_s\Lambda$ is identified with a space of $3\times 3$ symmetric matrices). Thus, the projectivised tangent bundle  of the hypersurface $M^5$ is equipped with a field of rational normal curves. Since the group of conformal automorphisms of the rational normal curve is isomorphic to $GL(2)$, this specifies on $M^5$ an irreducible $GL(2)$-structure. Alternatively, one can say that each tangent space to $M^5$ is canonically identified with a five-dimensional space of binary quartics; in this picture the rational normal curve corresponds to quartics with a quadruple root.  Similar structures were discussed in \cite {Dun, Bailey} under the name `paraconformal'.  

\medskip

{\bf Remark.} Another natural source of irreducible $GL(2)$-structures in dimension five is provided by 5-th order ordinary differential equations with the vanishing Doubrov invariants \cite{Doubrov, Dun}. Here $GL(2)$-structures arise on the moduli spaces of solutions to the corresponding ODEs. The  relation of these  structures to  equations of the  dispersionless Hirota type  is yet unclear. 

\medskip

Irreducible $GL(2)$-structures in dimension five possess a number of specific properties, allowing one to equip $M^5$ with  more `familiar' objects such as the conformal  cubic form and the conformal metric. These are defined as follows:

\noindent The {\bf conformal cubic form $C_{M^5}$} is  a restriction to $M^5$ of the  cubic form $C$ defined on $\Lambda$. The knowledge of this cubic form allows one to reconstruct the $GL(2)$-structure, indeed, the singular locus of the cubic $C_{M^5}$ is the rational normal curve introduced above. Conversely, given a rational normal curve, its bisecant variety is the cubic  $C_{M^5}=0$ \cite{Harris}, p.119. This specifies $C_{M^5}$ uniquely up to a conformal factor.

\noindent The {\bf conformal metric $Q_{M^5}$} is defined by the unique quadratic cone  in $T_sM^5$ containing the tangent variety of the rational normal curve, see \cite{Harris}, p.119,  for an alternative definition. The tangent variety of the rational normal curve is a complete intersection of  $Q_{M^5}$ and  $C_{M^5}$.

One can show that the conformal metric $Q_{M^5}$ and the cubic form $C_{M^5}$ satisfy a system of remarkable  relations which have appeared recently in \cite{Nur1} (notice that the signature of the quadratic form in \cite{Nur1} is Euclidean, while in our case the metric is pseudo-Euclidean of the signature $(3, 2)$).  Let $Q_{ij}$ and $C_{ijk}$ be the components of $Q_{M^5}$ and  $C_{M^5}$ in any coordinate system on $M^5$ (we point out that these objects are defined up to  conformal factors). One can show that the conformal factors can be specified in such a way that the following relations are satisfied:
\begin{equation}
\begin{array}{c}
C_{ijk}Q^{kj}=0 ~~~ ({\rm apolarity ~ condition}), \\
\ \\
C_{jkr}Q^{rs}C_{lns}+C_{ljr}Q^{rs}C_{kns}+C_{klr}Q^{rs}C_{jns}=Q_{jk}Q_{ln}+Q_{lj}Q_{kn}+
Q_{kl}Q_{jn};
\end{array}
\label{so3}
\end{equation}
these relations still allow simultaneous conformal reparametrizations of the form $C\to \varphi^3 C, \ Q\to \varphi^2 Q$.  
The explicit coordinate formulas for $Q_{M^5}$ and  $C_{M^5}$ are provided in Sect. 4.4. 

The normalised conformal metric $Q_{M^5}$ and  the cubic form $C_{M^5}$ satisfying the equations 
(\ref{so3}) carry  important geometric  information about the  hypersurface $M^5$. Their transformation laws suggest that one should consider the ratio $\sigma=(C_{M^5})^2/(Q_{M^5})^3$ which can be viewed as an analog of the Fubini projective element. 
By construction, this object is manifestly $Sp(6)$-invariant. We can formulate the following

\medskip

{\bf Conjecture.} {\it A generic five-dimensional hypersurface  of the Lagrangian Grassmanian $\Lambda$ is uniquely defined, modulo  $Sp(6)$-equivalence, by its  `symplectic element' $\sigma$.} 

\medskip

\subsection{Geometric interpretation of the integrability conditions: bisecant surfaces and trisecant $3$-folds}

As explained in Sect. 4.2, each projectivised tangent space $T_sM^5$ carries a rational normal curve $\gamma$ of degree four. Among the most natural geometric objects associated with a rational normal curve are its bisecant lines (two-parameter family) and trisecant planes (three-parameter family). 

\noindent {\bf Definition.} A {\it bisecant surface} is a two-dimensional submanifold $\Sigma^2\subset M^5$ whose projectivised tangent planes are bisecant lines. Similarly, a {\it trisecant $3$-fold} is a three-dimensional submanifold $\Sigma^3\subset M^5$ whose projectivised tangent spaces are trisecant planes.  To be more precise, we will consider {\it holonomic} trisecant $3$-folds which can be defined as follows. Notice first that each tangent space $T_s\Sigma^3$ carries three distinguished directions, namely, those corresponding to the three points of intersection of $PT_s\Sigma^3$ with $\gamma$.
These directions define a net on $\Sigma^3$, which we will call the {\it characteristic} net. We require  the characteristic net  to be  holonomic (that is,  a coordinate net). A similar net exists on any bisecant surface, however, the requirement of holonomicity is superfluous: in two dimensions, any net is automatically holonomic.

\begin{theorem} Bisecant surfaces and holonomic trisecant $3$-folds of a hypersurface $M^5$  correspond to two- and three-component hydrodynamic reductions of the associated equation of the dispersionless Hirota type. Moreover,

\noindent {\bf (i)} Each hypersurface $M^5$  possesses infinitely many bisecant surfaces parametrized  by two arbitrary functions of one variable; 

\noindent{\bf (ii)} A hypersuface $M^5$ corresponds to an integrable equation  if and only if it possesses infinitely many holonomic trisecant $3$-folds  parametrized by three arbitrary functions of one variable. Thus, the existence of holonomic trisecant $3$-folds is a geometric interpretation of the integrability property.
\end{theorem}

\medskip

\centerline {\bf Proof:}

We follow the notation of Sect. 2. Consider an equation of the dispersionless Hirota type represented in the form (\ref{evol}), and take its $n$-component reduction specified by  $a=a(R^{1},..., R^{n}),\ b=b(R^{1},...,R^{n}),\ c=c(R^{1},...,R^{n}),\ p=p(R^{1},...,R^{n}), \ q=q(R^{1},...,R^{n}), \ f=f(R^1,... R^n)$. The image of this solution in the Lagrangian Grassmannian is a submanifold $\Sigma^n\subset M^5$ (it would be sufficient for our purposes to restrict to   $n=2, 3$), represented by a symmetric matrix
\begin{equation}
U=
\left(\begin{array}{ccc}
a & b & p \\
b & c & q \\
p & q & f
\end{array}
\right)
\label{U}
\end{equation}
parametrized by $R^1, ..., R^n$. Using (\ref{bpqc}) and (\ref{disp}), one obtains the following expression for the derivative of $U$ with respect to $R^i$,
\begin{equation}
\partial_iU=
\left(\begin{array}{ccc}
1 & \mu^i & \lambda^i \\
\mu^i & (\mu^i)^2 & \lambda^i \mu^i \\
\lambda^i & \lambda^i \mu^i &( \lambda^i)^2
\end{array}
\right)\partial_ia;
\label{Ui}
\end{equation}
thus, $U_i$ is a matrix of rank one, so that the projectivization of $\partial_iU$ belongs to $\gamma$, and the coordinates $R^i$ provide a characteristic net on $\Sigma^n$. For $n=2$ we have a two-dimensional surface $\Sigma^2$ parametrized by $R^1, R^2$. Since both $\partial_1U$ and $\partial_2U$ have rank one,   the surface $\Sigma^2$ is bisecant. As follows from Sect. 2, any equation of the dispersionless Hirota type (not necessarily integrable) possesses infinitely many two-component reductions parametrized by two arbitrary functions of one variable. Thus, any hypersurface $M^5$ of the Lagrangian Grassmannian possesses an  infinity  of bisecant surfaces. This establishes the first part {(i)} of the theorem.

Similarly, any three-component reduction corresponds to a holonomic trisecant $3$-fold $\Sigma^3\subset M^5$ parametrized by $R^1, R^2, R^3$. Since an integrable equation possesses (by definition) infinitely many three-component reductions parametrized by three arbitrary functions of one variable, the corresponding hypersurface $M^5$ possesses an infinity of holonomic trisecant $3$-folds. This establishes the first part of (ii). To finish the proof one needs to show that, conversely,   bisecant surfaces (holonomic trisecant 3-folds) of $M^5$ corresponds to two-component (three-component) reductions of the associated equation. This can be demonstrated as follows. Let $\Sigma^2$ be a bisecant surface represented in the form (\ref{U}), referred to its characteristic net $R^1, R^2$. Thus,  $a, b, c, p, q, f$ are functions of $R^1, R^2$ such that the rank of $\partial_iU$ equals one, so that  one can introduce the parametrization (\ref{Ui}). The compatibility conditions for the equations $\partial_ib=\mu^i\partial_ia$ and $\partial_ip=\lambda^i\partial_ia$ imply
\begin{equation}
\partial_i\partial_ja=\frac{\partial_j\mu^i}{\mu^j-\mu^i}\partial_ia+\frac{\partial_i\mu^j}{\mu^i-\mu^j}\partial_ja ~~~ {\rm and} ~~~
\partial_i\partial_ja=\frac{\partial_j\lambda^i}{\lambda^j-\lambda^i}\partial_ia+\frac{\partial_i\lambda^j}{\lambda^i-\lambda^j}\partial_ja,
\label{x}
\end{equation} 
respectively. Using (\ref{x}), the consistency conditions for $\partial_ic=(\mu^i)^2\partial_ia$ and $\partial_if=(\lambda^i)^2\partial_ia$ simplify to 
$\partial_j\mu^i\partial_ia+\partial_i\mu^j\partial_ja=0$ and $\partial_j\lambda^i\partial_ia+\partial_i\lambda^j\partial_ja=0$, compare with (\ref{Comp}). Substituting these relations into (\ref{x}) one obtains
$$
\partial_i\partial_ja=2\frac{\partial_j\mu^i}{\mu^j-\mu^i}\partial_ia ~~~ {\rm and} ~~~
\partial_i\partial_ja=2\frac{\partial_j\lambda^i}{\lambda^j-\lambda^i}\partial_ia,
$$
which implies the commutativity condition (\ref{Com}). Finally, the last compatibility condition,
$\partial_iq=\mu^i\lambda^i\partial_ia$,  holds identically. The dispersion relation (\ref{disp}), which connects $\mu^i$ and $\lambda^i$, follows from the relation $\partial_if=(\lambda^i)^2\partial_ia$
after the substitution $f=f(a, b, c, p, q)$. Ultimately, we have recovered all equations defining hydrodynamic reductions, see Sect. 2. Thus, bisecant surfaces correspond to two-component reductions  (for trisecant $3$-folds, considerations are exactly the same).

We have characterized the integrability as the existence of `sufficiently many' special submanifolds. This is very reminiscent  of the characterization of self-dual $4$-manifolds in terms of the existence of $\alpha$-surfaces \cite{Atiyah}, or a similar characterization of  semi-integrable almost Grassmann structures \cite{Akivis2}.

\subsection{Computational formulas}

Here we provide the explicit coordinate formulae for the conformal metric and the conformal cubic form as introduced in Sect. 4.2.

The Lagrangian Grassmannian $\Lambda$ is (locally) parametrized  by $3\times 3$ symmetric matrices (\ref{U}), so that the cubic form $C$ is given by
$$
C= \det dU= dadcdf-da(dq)^2-(db)^2df+2dbdpdq-dc(dp)^2.
$$
Each tangent space to the Lagrangian Grassmannian (which can also be identified with $3\times 3$ symmetric matrices)  carries a cubic cone $C=0$ consisting of rank two symmetric matrices. Its singular locus coincides with rank one matrices, and can be parametrized in the form
\begin{equation}
\left(\begin{array}{ccc}
1 & \mu & \lambda \\
\mu & \mu^2 & \lambda \mu \\
\lambda & \lambda \mu & \lambda^2
\end{array}
\right)
\label{Ver}
\end{equation}
where $\lambda$ and $\mu$ are independent parameters. The restriction $C_{M^5}$ of the cubic form $C$ to the hypersurface $M^5$,  defined by the equation $f=f(a, b, c, p, q)$, is obtained by setting 
\begin{equation}
df=f_ada+f_bdb+f_cdc+f_pdp+f_qdq.
\label{df}
\end{equation}
This results in
\begin{equation}
\begin{array}{c}
C_{M^5}=f_a(da)^2dc+f_bdadbdc+f_pdadcdp+f_cda(dc)^2+f_qdadcdq-da(dq)^2-\\
\ \\
f_ada(db)^2-f_b(db)^3-f_p(db)^2dp-f_c(db)^2dc-f_q(db)^2dq+2dbdpdq-dc(dp)^2,
\end{array}
\label{C}
\end{equation}
here $a, b, c, p, q$ are  local coordinates on $M^5$. The intersection of the singular locus (\ref{Ver}) with the tangent space to $M^5$ is a rational normal curve with the  affine parametrization
$$
\gamma=(1, \ \mu, \ \mu^2, \ \lambda, \ \mu \lambda);
$$
notice that $\lambda$ and $\mu$ are no longer independent and satisfy, by virtue of (\ref{df}),  the quadratic relation
$$
\lambda^2=f_{a}+f_{b}\mu+f_{c}\mu^2+f_{p}\lambda+f_{q}\lambda \mu,
$$
which coincides with the dispersion relation (\ref{disp}). The tangent variety of  $\gamma$ is defined as $\gamma+t \gamma'$ where prime denotes differentiation with respect to $\mu$, and $\lambda$ is viewed as a function of $\mu$ specified by the dispersion relation. Explicitly, the tangent variety is given by
$$
(1, \ \mu+t, \ \mu^2+2t\mu, \ \lambda+t\lambda', \ \mu \lambda+t(\lambda+\mu \lambda');
$$
here $\lambda'$  is obtained by implicitly differentiating the dispersion relation with respect to $\mu$,
$$
\lambda'=(f_b+2f_c\mu+f_q\lambda)/(2\lambda-f_q\mu-f_p).
$$
One can verify that, up to a conformal factor,  there exists a unique  quadratic form $Q_{M^5}$ vanishing on the tangent variety of $\gamma$,
\begin{equation}
\begin{array}{c}
Q_{M^5}=(4f_a^2+f_af_p^2)da^2+(8f_af_b+f_bf_p^2+2f_af_pf_q)dadb+(4f_af_p+f_p^3)dadp\\
\ \\
+(f_b^2+4f_af_c+f_cf_p^2+f_af_q^2)dadc+(2f_bf_p+f_p^2f_q)dadq \\
\ \\
+(3f_b^2+4f_af_c+2f_bf_pf_q)db^2+(2f_bf_p+4f_af_q+2f_p^2f_q)dbdp \\
\ \\
+(8f_bf_c+2f_cf_pf_q+f_bf_q^2)dbdc+(4f_cf_p+2f_bf_q+2f_pf_q^2)dbdq\\
\ \\
-(4f_a+f_p^2)dp^2+(2f_bf_q+f_pf_q^2)dpdc-(4f_b+2f_pf_q)dpdq\\
\ \\
+(4f_c^2+f_cf_q^2)dc^2+(4f_cf_q+f_q^3)dcdq-(4f_c+f_q^2)dq^2;
\end{array}
\label{Q}
\end{equation}
(one needs to use the dispersion relation in order to verify the vanishing of $Q_{M^5}$ on the tangent variety). Representing $Q_{M^5}$ as a $5\times 5$ symmetric matrix and calculating its determinant one obtains
$$
\det Q_{M^5}=3(f_b^2+f_bf_pf_q-f_af_q^2-f_cf_p^2-4f_af_c)^4=3\triangle^4;
$$
we point out that  $\triangle = 0$ if and only if the corresponding dispersion relation (\ref{disp}) defines a reducible conic. Thus, under the assumption of the irreducibility, the conformal metric $Q_{M^5}$ is non-degenerate.  The cubic form (\ref{C}) and the conformal metric (\ref{Q}) satisfy the relations
$$
\begin{array}{c}
C_{ijk}Q^{kj}=0, \\
\ \\
C_{jkr}Q^{rs}C_{lns}+C_{ljr}Q^{rs}C_{kns}+C_{klr}Q^{rs}C_{jns}=\frac{1}{27\triangle^2}(Q_{jk}Q_{ln}+Q_{lj}Q_{kn}+
Q_{kl}Q_{jn}).
\end{array}
$$
Thus, the normalization $C\to 3\sqrt3 \ \triangle C$ results in the identities (\ref{so3}).

\subsection{General aspects of irreducible $Gl(2)$-structures}

 Let $P^4$ be a projective space with coordinates $(x^0 : x^1 : x^2 : x^3 : x^4)$.  A rational normal curve $\gamma$ of degree four can be parametrized in the form $\gamma=(1 : t : t^2 : t^3 : t^4)$. Its bisecant variety is a cubic hypersurface defined by the equation
$$
 \det \left(
\begin{array}{ccc}
x^0 & x^1 & x^2 \\
x^1 & x^2 & x^3 \\
x^2 & x^3 & x^4 \\
\end{array}
\right)=0;
$$
notice that $\gamma$ can be recovered as a singular locus of its bisecant variety. The tangent variety of $\gamma$ is contained in a unique quadric hypersurface 
$$
 x^0x^4-4x^1x^3+3(x^2)^2=0.
$$

The `curved' analog of this picture is the following. Let $M^5$ be a $5$-dimensional manifold with a field of rational normal curves of degree four specified in the projectivization of each tangent space. One can choose a frame of $1$-forms $\omega^1, \omega^2, \omega^3, \omega^4, \omega^5$ such that the equation of  the bisecant variety of $\gamma$ takes the canonical  form $C=0$ where
$$
C= \det \left(
\begin{array}{ccc}
\omega^1 & \omega^2 & \omega^3 \\
\omega^2 & \omega^3 & \omega^4 \\
\omega^3 & \omega^4 & \omega^5 \\
\end{array}
\right);
$$
such frame is defined up to a natural $GL(2)$-equivalence,  leading to the following structure equations:
\begin{eqnarray}
&&d\omega^1=\Omega \wedge \omega^1+2\varphi \wedge \omega^1+4\eta \wedge \omega^2+...,\notag 
\\
&&d\omega^2=\Omega \wedge \omega^2+\varphi \wedge \omega^2+\psi \wedge \omega^1+3\eta\wedge \omega^3+...,\notag \\
&&d\omega^3=\Omega \wedge \omega^3+2\psi \wedge \omega^2+2\eta \wedge \omega^4+...,
\label{omega} \\
&&d\omega^4=\Omega \wedge \omega^4-\varphi \wedge \omega^4+3\psi \wedge \omega^3+ \eta\wedge \omega^5+...,\notag \\
&&d\omega^5=\Omega \wedge \omega^5-2\varphi \wedge \omega^5+4\psi \wedge \omega^4+....\notag
\end{eqnarray}
Here dots denote linear combinations of the terms $\omega^i\wedge \omega^j$ which are responsible for the `torsion'. The structure equations for the secondary forms $\Omega, \varphi, \psi, \eta$,  restricted to the fiber $\omega^i=0$, take the standard  $GL(2)$-form:
$$
d\Omega=0, ~~~ d\eta=\varphi \wedge \eta, ~~~ d\varphi=2\eta \wedge \psi, ~~~ d\psi=-\varphi \wedge \psi.
$$
Introducing ${\bf \omega}= (\omega^1, \omega^2, \omega^3, \omega^4, \omega^5)^t $, one can represent (\ref{omega}) in a compact form
\begin{equation}
d{\bf \omega}=(\Omega E+\psi E_1+\varphi H+ \eta E_2)\wedge {\bf \omega}+T({\bf \omega}, {\bf \omega}),
\label{compact}
\end{equation}
where $T$ represents the torsion, $E$ is the identity matrix, and $E_1, H, E_2$ define the $5$-dimensional irreducible representation of $SL(2)$:
$$
E_1=\left(
\begin{array}{ccccc}
0&0&0&0&0 \\
1&0&0&0&0 \\
0&2&0&0&0 \\
0&0&3&0&0 \\
0&0&0&4&0 
\end{array}
\right), ~~~
H=\left(
\begin{array}{ccccc}
2&0&0&0&0 \\
0&1&0&0&0 \\
0&0&0&0&0 \\
0&0&0&-1&0 \\
0&0&0&0&-2 
\end{array}
\right), ~~~
E_1=\left(
\begin{array}{ccccc}
0&4&0&0&0 \\
0&0&3&0&0 \\
0&0&0&2&0 \\
0&0&0&0&1 \\
0&0&0&0&0 
\end{array}
\right).
$$
We point out that bisecant surfaces and trisecant three-folds make perfect sense for abstract $GL(2)$ structures. 
Recall that if the $GL(2)$ structure  comes from a hypersurface of the Lagrangian Grassmanian, it automatically possesses infinitely many bisecant surfaces parametrized by two arbitrary functions of one variable. For abstract $GL(2)$ structures, the existence of bisecant surfaces imposes strong restrictions on the torsion $T$, which can be represented in a simple geometric form:

\begin{theorem} An abstract $GL(2)$ structure possesses infinitely many bisecant surfaces  if and only if,  for any two vectors $X$ and $Y$
which belong to the rational normal curve, one has $C(X, Y, T(X, Y))=0$.
\end{theorem}

\centerline{\bf Proof:}
Let us  look for bisecant surfaces in the form
$
\omega=\alpha\  {\bf a}\  dR^1+\beta\  {\bf b}\  dR^2
$
where $R^1, R^2$ are the characteristic coordinates,  ${\bf a}$ and ${\bf b}$ are (column)  vectors which belong to the rational normal curve $\gamma$, that is, 
$$
{\bf a}=(1 : a :  a^2 : a^3 : a^4)^t, ~~~ {\bf b}=(1 : b :  b^2 : b^3 : b^4)^t,
$$
and $a, \ b, \ \alpha, \ \beta$ are certain functions of $R^1, R^2$. Substituting this ansatz into (\ref{compact}) and using the identities
\begin{eqnarray}
&&E_1{\bf a}={\bf a}',\quad
E_1{\bf b}={\bf b}',\notag \\
&&H{\bf a}=2{\bf a}-a{\bf a}',\quad
H{\bf b}=2{\bf b}-b{\bf b}',\notag \\
&&E_2{\bf a}=4a{\bf a}-a^2{\bf a}',\quad
E_2{\bf b}=4b{\bf b}-b^2{\bf b}',\notag
\end{eqnarray}
where
$$
{\bf a}'=(0 : 1 :  2a : 3a^2 : 4a^3)^t, ~~~ {\bf b}=(0 : 1 :  2b : 3b^2 : 4b^3)^t,
$$
one arrives at 
\begin{equation}
\begin{array}{c}
{\bf a}(d\alpha-\alpha \Omega-2\alpha \varphi-4a\alpha \eta)\wedge dR^1+
{\bf b}(d\beta-\beta\Omega-2\beta \varphi-4b\beta \eta)\wedge dR^2+ \\
\ \\
{\bf a}'(\alpha da-\alpha \psi+\alpha a \varphi+\alpha a^2\eta)\wedge dR^1+
{\bf b}'(\beta db-\beta \psi+\beta b \varphi+\beta b^2\eta)\wedge dR^2= \\
\ \\
2\alpha \beta \ T({\bf a}, {\bf b}) \ dR^1\wedge dR^2.
\end{array}
\label{+}
\end{equation}
Thus, for the consistency of these equations one has to require
$$
T({\bf a}, {\bf b})\in {\rm span} \ \{{\bf a}, {\bf b}, {\bf a}', {\bf b}'\},
$$
which is equivalent to the statement of the theorem. Once this condition is satisfied, one obtains a system of four first order PDEs for $a, b, \alpha, \beta$ by collecting similar terms in (\ref{+}). Up to reparametrizations $R^i\to \varphi^i(R^i)$, the general solution of this system  depends on two arbitrary functions of one variable.

Theorem 4 establishes  necessary conditions for the realizability of a $GL(2)$ structure on a hypersurface of the Lagrangian Grassmannian. Further development of the general theory of abstract $GL(2)$ structures  is beyond the scope of this paper, and will be addressed elsewhere.

\section{Integrability in more than three dimensions} 

The integrability of three-dimensional equations of the form (\ref{main}) was defined as the existence of n-component hydrodynamic reductions (\ref{R}) parametrized by $n$ arbitrary functions of a single variable. As demonstrated in \cite{Fer22, Fer3}, this approach readily generalizes to any dimension: a $d$-dimensional PDE
\begin{equation}
F(u_{x^ix^j})=0
\label{4}
\end{equation}
for a function $u$ of $d$ independent variables $x^1, ..., x^d$ is said to be integrable if it possesses infinitely many $n$-component hydrodynamic reductions  parametrized by $(d-2)n$ arbitrary functions of a single variable. In this case  equations (\ref{R}) are replaced  by $d-1$ commuting $(1+1)$-dimensional systems of hydrodynamic type. Among the known four-dimensional integrable examples one should primarily mention the `first heavenly equation', 
$$
u_{xy}u_{zt}-u_{xt}u_{zy}=1,
$$
as well as  its equivalent forms,
$$
u_{tx}+u_{zy}+u_{xx}u_{yy}-u^2_{xy}=0
$$
and
$$
u_{tt}=u_{xy}u_{zt}-u_{xt}u_{zy},
$$
known as the `second heavenly' and the `Grant equation', respectively \cite{Plebanski, Grant}. It was demonstrated in \cite{Fer22, Fer3} that $n$-component reductions of these equations are parametrized by $2n$ arbitrary functions of a single variable. 

An interesting six-dimensional integrable generalization of the heavenly equation,
$$
u_{t\tilde t}+u_{z\tilde z}+u_{tx}u_{zy}-u_{ty}u_{zx}=0,
$$
arises in the context of $sdiff(\Sigma^2)$  self-dual Yang-Mills equations \cite{Przanovski}.  Its $n$-component  reductions are parametrized by $4n$ arbitrary functions of a single variable \cite{Fer3}. Notice that all these examples belong to the class of symplectic Monge-Amp\'ere equations as introduced in Sect. 3.

Although the general problem of classification of multi-dimensional integrable equations can be approached in a similar way, the method of Sect. 2 leads to quite complicated analysis. One way to 
bypass lengthy calculations is based on the following simple idea: suppose we want to classify four-dimensional integrable equations of the form (\ref{4}) for a function $u(x, y, z, t)$. Let us look for travelling wave solutions in the form
$$
u(X, Y, Z)=u(x+\alpha t, \ y+\beta t, \ z+\gamma t).
$$
The substitution of this ansatz into (\ref{4}) leads to a three-dimensional equation which must be integrable for {\it any} values of constants $\alpha, \beta, \gamma$. Since,  in  three dimensions, the integrability conditions are explicitly known, this provides strong restrictions on the original function $F$, which are {\it necessary} for the integrability. The philosophy of this approach is well familiar from the soliton theory: symmetry reductions of integrable systems must be themselves integrable.

Thus, for the first heavenly equation, traveling wave solutions are governed by
$$
\alpha(u_{XY}u_{XZ}-u_{XX}u_{YZ})+\gamma(u_{XY}u_{ZZ}-u_{XZ}u_{YZ})=1,
$$
which is a three-dimensional symplectic Monge-Ampere equation. One can show that it is indeed integrable (in fact, linearizable)
 for any values of constants.

We hope to return to the multi-dimensional case in subsequent publications.

\section{$Sp(6)$-action on the moduli space  of integrable equations}

In this section we investigate the action of the equivalence group $Sp(6)$ on the moduli space of integrable equations of the dispersionless Hirota type (recall that both $Sp(6)$ and the moduli space have coinciding dimensions equal to $21$). Our main result states that this action has an open orbit.

For a $3\times 3$ symmetric matrix $U=[u_{ij}] $, the Lie algebra of the group of conformal automorphisms of the cubic form $det \ dU$  is spanned by  21 vector fields which generate the Lie algebra of the symplectic Lie group $Sp(6)$:

\begin{eqnarray*}
&&{\bf X_{11}}=\frac{\partial}{\partial{u_{11}}} ~,~ {\bf X_{12}}=\frac{\partial}{\partial{u_{12}}} ~,~ {\bf X_{13}}=\frac{\partial}{\partial{u_{13}}} ~,~ {\bf X_{22}}=\frac{\partial}{\partial{u_{22}}}~,~ {\bf X_{23}}=\frac{\partial}{\partial{u_{23}}} ~,~ {\bf X_{33}}=\frac{\partial}{\partial{u_{33}}},
\end{eqnarray*}

\begin{eqnarray*}
&&{\bf J_1}=2 u_{11}\frac{\partial}{\partial{u_{11}}}+u_{12}\frac{\partial}{\partial{u_{12}}}+u_{13}\frac{\partial}{\partial{u_{13}}},  \notag \\
&&{\bf J_2}=2 u_{22}\frac{\partial}{\partial{u_{22}}}+u_{21}\frac{\partial}{\partial{u_{21}}}+u_{23}\frac{\partial}{\partial{u_{23}}},  \notag \\
&&{\bf J_3}=2 u_{33}\frac{\partial}{\partial{u_{33}}}+u_{31}\frac{\partial}{\partial{u_{31}}}+u_{32}\frac{\partial}{\partial{u_{32}}},  
\end{eqnarray*}

\begin{eqnarray*}
&&{\bf L_{12}}=2 u_{12}\frac{\partial}{\partial{u_{11}}}+u_{22}\frac{\partial}{\partial{u_{12}}}+u_{23}\frac{\partial}{\partial{u_{13}}},  \notag \\
&&{\bf L_{13}}=2 u_{13}\frac{\partial}{\partial{u_{11}}}+u_{33}\frac{\partial}{\partial{u_{13}}}+u_{32}\frac{\partial}{\partial{u_{12}}},  \notag \\
&&{\bf L_{21}}=2 u_{21}\frac{\partial}{\partial{u_{22}}}+u_{11}\frac{\partial}{\partial{u_{21}}}+u_{13}\frac{\partial}{\partial{u_{23}}},  \notag \\
&&{\bf L_{23}}=2 u_{23}\frac{\partial}{\partial{u_{22}}}+u_{33}\frac{\partial}{\partial{u_{23}}}+u_{31}\frac{\partial}{\partial{u_{21}}},  \notag \\
&&{\bf L_{31}}=2 u_{31}\frac{\partial}{\partial{u_{33}}}+u_{11}\frac{\partial}{\partial{u_{31}}}+u_{12}\frac{\partial}{\partial{u_{23}}},  \notag \\
&&{\bf L_{32}}=2 u_{32}\frac{\partial}{\partial{u_{33}}}+u_{22}\frac{\partial}{\partial{u_{32}}}+u_{21}\frac{\partial}{\partial{u_{13}}},  
\end{eqnarray*}

\begin{eqnarray*}
&&{\bf H_1}={u_{11}}^2\frac{\partial}{\partial{u_{11}}}+u_{11}u_{12}\frac{\partial}{\partial{u_{12}}}+u_{11}u_{13}\frac{\partial}{\partial{u_{13}}}+{u_{12}}^2\frac{\partial}{\partial{u_{22}}}+u_{12}u_{13}\frac{\partial}{\partial{u_{23}}}+{u_{13}}^2\frac{\partial}{\partial{u_{33}}},  \notag \\
&&{\bf H_2}={u_{22}}^2\frac{\partial}{\partial{u_{22}}}+u_{22}u_{21}\frac{\partial}{\partial{u_{21}}}+u_{22}u_{23}\frac{\partial}{\partial{u_{23}}}+{u_{21}}^2\frac{\partial}{\partial{u_{11}}}+u_{21}u_{23}\frac{\partial}{\partial{u_{13}}}+{u_{23}}^2\frac{\partial}{\partial{u_{33}}},  \notag \\
&&{\bf H_3}={u_{33}}^2\frac{\partial}{\partial{u_{33}}}+u_{33}u_{31}\frac{\partial}{\partial{u_{31}}}+u_{33}u_{32}\frac{\partial}{\partial{u_{32}}}+{u_{31}}^2\frac{\partial}{\partial{u_{11}}}+u_{31}u_{32}\frac{\partial}{\partial{u_{21}}}+{u_{32}}^2\frac{\partial}{\partial{u_{22}}}, 
\end{eqnarray*}

\begin{eqnarray*}
&&{\bf P_1}=2 u_{12}u_{13}\frac{\partial}{\partial{u_{11}}}+(u_{12}u_{23}+u_{13}u_{22})\frac{\partial}{\partial{u_{12}}}+(u_{12}u_{33}+u_{13}u_{23})\frac{\partial}{\partial{u_{13}}}+2 u_{22}u_{23}\frac{\partial}{\partial{u_{22}}}  \\
&&+(u_{22}u_{33}+{u_{23}}^2)\frac{\partial}{\partial{u_{23}}}+2 u_{23}u_{33}\frac{\partial}{\partial{u_{33}}}, \notag \\
&&{\bf P_2}=2 u_{21}u_{23}\frac{\partial}{\partial{u_{22}}}+(u_{21}u_{13}+u_{23}u_{11})\frac{\partial}{\partial{u_{21}}}+(u_{21}u_{33}+u_{23}u_{13})\frac{\partial}{\partial{u_{23}}}+2 u_{11}u_{13}\frac{\partial}{\partial{u_{11}}} \\ 
&&+(u_{11}u_{33}+{u_{13}}^2)\frac{\partial}{\partial{u_{13}}}+2 u_{13}u_{33}\frac{\partial}{\partial{u_{33}}},  \notag \\
&&{\bf P_3}=2 u_{31}u_{32}\frac{\partial}{\partial{u_{33}}}+(u_{31}u_{12}+u_{32}u_{11})\frac{\partial}{\partial{u_{31}}}+(u_{31}u_{22}+u_{32}u_{12})\frac{\partial}{\partial{u_{32}}}+2 u_{11}u_{12}\frac{\partial}{\partial{u_{11}}}\\
&&+(u_{11}u_{22}+{u_{12}}^2)\frac{\partial}{\partial{u_{12}}}+2 u_{12}u_{22}\frac{\partial}{\partial{u_{22}}}.\notag
\end{eqnarray*}
Given a PDE of the form  $F(u_{11}, u_{12}, u_{13}, u_{22}, u_{23}, u_{33})=0$, we will look for its infinitesimal symmetries by solving the determining equation $L_X F\vert_{F=0}=0$, where $X$ is a linear combination of the $21$ vector fields presented above (notice that the answer may not coincide with the full algebra of Lie-point symmetries: we consider only those symmetries which belong to the equivalence group $Sp(6)$). Below we list particular examples of integrable equations which possess symmetry algebras of different dimensions (it is worth noting that any two equations with different symmetry algebras are automatically non-equivalent).

\noindent {\bf Example 1.} The 2-dimensional  linear wave equation, $u_{11}+u_{22}-u_{33}=0$, possesses nine infinitesimal symmetries:
$$
{\bf X}_{12}, ~~ {\bf X}_{13}, ~~ {\bf X}_{23}, ~~ {\bf X}_{11}+ {\bf X}_{33}, ~~ {\bf X}_{22}+ {\bf X}_{33}, ~~ 
{\bf J}_{1}+{\bf J}_{2}+{\bf J}_{3}, ~~ {\bf L}_{12}-{\bf L}_{21}, ~~  {\bf L}_{13}+     {\bf L}_{31}, ~~ {\bf L}_{23}+{\bf L}_{32}. 
$$
One can prove that the existence of nine infinitesimal symmetries is necessary and sufficient for the linearizability  of a  general equation of the form  (\ref{main}),  \cite{Nur}.

\noindent {\bf Example 2.} The dKP equation, $u_{22}-u_{13}+\frac{1}{2}u^2_{11}=0$, possesses seven  infinitesimal symmetries:
$$
{\bf X}_{12}, ~~ {\bf X}_{23}, ~~ {\bf X}_{33}, ~~ {\bf X}_{13}+ {\bf X}_{22}, ~~ {\bf X}_{11}+ {\bf L}_{31}, ~~ 
{\bf J}_{1}+2{\bf J}_{2}+3{\bf J}_{3}, ~~ 2{\bf L}_{32}+{\bf L}_{21}. 
$$
It is likely that there exist no integrable equations with eight symmetries.

\noindent {\bf Example 3.} The Boyer-Finley equation, $u_{11}+u_{22}-e^{u_{33}}=0$, possesses six   infinitesimal symmetries:
$$
{\bf X}_{12}, ~~ {\bf X}_{13}, ~~ {\bf X}_{23}, ~~ {\bf X}_{11}- {\bf X}_{22}, ~~ 
{\bf J}_{1}+{\bf J}_{2}+2{\bf X}_{33}, ~~ {\bf L}_{12}-{\bf L}_{21}. 
$$

\noindent {\bf Example 4.} The degeneration of the dispersionless Hirota equation, $u_{12}-u_{13}-e^{u_{23}}=0$,  possesses five  infinitesimal symmetries:
$$
{\bf X}_{11}, ~~ {\bf X}_{22}, ~~ {\bf X}_{33}, ~~ {\bf J}_{1}+ {\bf X}_{23}, ~~ {\bf X}_{12}+ {\bf X}_{13}. 
$$

\noindent {\bf Example 5.} The   equation  $e^{u_{11}}+e^{u_{22}}-e^{u_{33}}=0$ possesses four  infinitesimal symmetries:
$$
{\bf X}_{12}, ~~ {\bf X}_{13}, ~~ {\bf X}_{23}, ~~ {\bf X}_{11}+ {\bf X}_{22}+ {\bf X}_{33}. 
$$

\noindent {\bf Example 6.} The dispersionless Hirota-type equation for the BKP hierarchy,
$e^{ u_{13}}+e^{ u_{13}}+e^{ u_{23}}=e^{ u_{13}+ u_{13}+ u_{23}}$, possesses three infinitesimal symmetries:
$$
{\bf X}_{11}, ~~ {\bf X}_{22}, ~~ {\bf X}_{33}. 
$$

\noindent {\bf Remark.} We point out that the existence of  `many'  symmetries is  not related to the integrability: for instance,  both equations (\ref{MA}), which are not integrable, possess $8$-dimensional symmetry algebras isomorphic to $SL(3, R)$ and $SU(3, R)$, respectively \cite{Banos}. Thus, the equation  $Hess~u=1$  possesses eight  infinitesimal symmetries:
$$
{\bf L}_{12}, ~~ {\bf L}_{13}, ~~ {\bf L}_{21}, ~~ {\bf L}_{23}, ~~ {\bf L}_{31}, ~~ {\bf L}_{32}, ~~ 
{\bf J}_{1}-{\bf J}_{2}, ~~ {\bf J}_{1}-{\bf J}_{3}.
$$
The equation  $Hess~u=u_{11}+u_{22}+u_{33}$  also possesses eight  infinitesimal symmetries:
$$
{\bf X}_{22}-{\bf X}_{11}, ~~ {\bf X}_{33}-{\bf X}_{11}, ~~ {\bf P}_{1}+{\bf X}_{23}~~{\bf P}_{2}+{\bf X}_{13}, ~~{\bf P}_{3}+{\bf X}_{12}, ~~{\bf L}_{12}+{\bf L}_{21}, ~~{\bf L}_{13}+{\bf L}_{31}, ~~{\bf L}_{23}+{\bf L}_{32}. 
$$

The main result of this section is the following

\begin{theorem}

The action of the equivalence group $Sp(6)$ on the moduli space of integrable equations of the dispersionless Hirota type has an open orbit.

\end{theorem}

This fact is, in a sense, surprising: it establishes the existence of a `universal' equation with no symmetries, which generates an open part of the moduli space under the action of $Sp(6)$. In particular, one should be able to obtain all equations with non-trivial symmetries by taking appropriate degenerations of this universal  equation.

\medskip

\centerline {\bf {Proof of Theorem 2}} 

\medskip

The main idea of the proof is to prolong the $21$ infinitesimal generators ${\bf X}_{11} - {\bf P}_3$ to the moduli space  of solutions of the involutive system (\ref{*}). We point out that, since third order derivatives of $f$ are explicitly known, this moduli space can be identified with the values of $f$ and its partial derivatives $f_i$, $f_{ij}$ up to second order ($21$ parameters altogether). 
The prolongation can be calculated as follows:

\noindent (1) Following the standard notation adopted in the symmetry analysis of differential equations
\cite{Ibragimov, Olver}, we introduce the variables
$$
x^1=u_{11}, ~~  x^2=u_{12}, ~~ x^3=u_{13}, ~~ x^4=u_{22}, ~~ x^5=u_{23}, ~~ u=u_{33},
$$
and rewrite the above $21$ generators in the form
$$
 \xi^i\frac{\partial}{\partial x^i}+\eta \frac{\partial}{\partial u};
$$
here $\xi^i$ and $\eta$ are certain functions of $x$ and $u$. In this notation, a dispersionless Hirota-type equation is represented in the form $u=u(x^1, ..., x^5)$ (the function $u$ is denoted by $f$ in Sect. 2). 

\noindent (2) Prolong   infinitesimal generators to the second order jet space with coordinates $u, u_i, u_{ij}$, 
$$
 \xi^i\frac{\partial}{\partial x^i}+\eta \frac{\partial}{\partial u}+
 \zeta_i\frac{\partial}{\partial u_i}+
 \zeta_{ij}\frac{\partial}{\partial u_{ij}},
$$
where $\zeta_i$ and $\zeta_{ij}$ are calculated according to the standard prolongation formulae
$$
\zeta_i=D_i(\eta)-u_kD_i(\xi^k), ~~~ \zeta_{ij}=D_j{\zeta_i}-u_{ik}D_j(\xi^k);
$$
here $D_i$ are the operators of total differentiation.

\noindent (3) To eliminate the $\frac{\partial}{\partial x^i}$-terms, 
subtract  the linear combination of total derivatives
$ \xi^i D_i$ from the prolonged operators. It is sufficient to keep only the following terms in $D_i$:
$$
D_i=\frac{\partial}{\partial x^i}+u_i \frac{\partial}{\partial u}+u_{ij}\frac{\partial}{\partial u_j}+u_{ijk}\frac{\partial}{\partial u_{jk}};
$$
notice that, since $u_{ijk}$ are explicit functions of lower order derivatives, the resulting operators will be 
well-defined on the $21$-dimensional space with coordinates $u, u_i, u_{ij}$. Although these operators  will depend on the variables $x^i$ as parameters (indeed, the isomorphism of the moduli space with the space $u, u_i, u_{ij}$ depends on the choice of a point in the $x$-space), all algebraic properties of these operators will be $x$-independent. 

\noindent (4) Finally, the dimension of the maximal $Sp(6)$-orbit  equals the rank of the $21\times 21$ matrix  of  coefficients of these operators. It remains to point out that this rank equals $21$ for any `random'  choice of numerical  values for $x^i, u, u_i, u_{ij}$ (however, it equals $21-r$ for any example with $r$ symmetries).

\section{Concluding remarks}

We would like to formulate a list of natural questions which are left beyond the scope of this paper.

\noindent --- It was demonstrated that the integrability of  equations of the form (\ref{main}) can be interpreted as the existence of `sufficiently many'  trisecant submanifolds of the corresponding hypersurfaces in the Lagrangian Grassmannian. It would be desirable to have a tensor characterization of this condition in terms of the metric $Q_{ij}$ and the cubic form $C_{ijk}$ introduced in Sect. 4.2. Presumably, these conditions will involve the Weyl tensor of the conformal metric $Q_{ij}$ (which is not conformally flat in general).
Alternatively, one may try to express the integrability conditions using the structure equations of the corresponding $GL(2)$-structure.

\noindent ---  We have shown that any hypersurface of the Lagrangian Grassmannian possesses an intrinsic $GL(2)$-structure. In would be desirable to develop a general theory of hypersurfaces of the Lagrangian Grassmannian, as well as the  theory of abstract $GL(2)$-structures, and to establish the embedding theorem. Notice that the theory of curves in the Lagrangian Grassmannian received some attention, see \cite{Ovs, Befa, Z}.

\noindent --- One of our main observations is the existence of a universal `master-equation'  which generates all other integrable equations of the dispersionless Hirota type via various (singular) limits. It seems to be very important to understand the reduction procedure geometrically, and to obtain  the master-equation explicitly in terms of the appropriate special functions (one shouldn't expect  simple formulas since even particular examples contain transcendental functions such as the elliptic theta-functions and solutions to the Chazy equation). 

To some extent, the picture resembles the situation with the classification of (non-holonomic) cyclids of Dupin in $E^4$: as demonstrated in \cite{Pinkall}, up to the action of the Lie sphere group $SO(5, 2)$, there exists a unique representative, namely,   the stereographic projection of the isoparametric hypersurface of Cartan. Coinciding dimensions of the equivalence groups, dim $Sp(6)$=dim $SO(5, 2)$=21, add to this similarity. Notice, however, that any nonholonomic Dupin hypersurface in $E^4$ is Lie-homogeneous, while our generic example possesses no continuous $Sp(6)$-symmetries.

\section*{Acknowledgements}
We thank Dmitry Alekseevsky, Robert Bryant, Fran Burstall, David Calderbank, Boris Doubrov, Nigel Hitchin,  Joseph Landsberg,  Pavel Nurowski, Maxim Pavlov and Volodya Roubtsov for clarifying discussions. The research of EVF was partially supported by the EPSRC grant EP/D036178/1, the European Union through the FP6 Marie Curie RTN project {\em ENIGMA} (contract number MRTN-CT-2004-5652), and the ESF programme MISGAM.

\end{document}